\documentclass[12pt,a4paper]{article}
\usepackage{amsthm}
\usepackage{fleqn}
\usepackage{array}
\usepackage[latin1]{inputenc}
\usepackage{amssymb}
\usepackage{amsfonts}
\usepackage{amsmath}

\begin{document}

\newtheorem{thm}{Theorem}[section]
\newtheorem{prop}{Proposition}[section]
\newtheorem{lem}{Lemma}[section]
\newtheorem{cor}{Corollary}[section]

\begin{center}
{\large {\bf Differential geometry of spatial curves for gauges}}

\vspace{5mm}

Vitor Balestro, Horst Martini and Makoto Sakaki
\end{center}

\vspace{5mm}

{\bf Abstract.} We derive Frenet-type results and invariants of spatial curves immersed
in $3$-dimensional generalized Minkowski spaces, i.e., in linear spaces which
satisfy all axioms of finite dimensional real Banach spaces except for
the symmetry axiom. Further on, we characterize cylindrical helices and
rectifying curves in such spaces, and the computation of invariants is
discussed, too. Finally, we study how translations of unit spheres
influence invariants of spatial curves.

\vspace{3mm}

{\bf Mathematics Subject Classification (2010).} 46B20, 52A15, 52A21, 53A04, 53A35

\vspace{3mm}

{\bf Keywords.} gauges, generalized Minkowski spaces, helix, invariants, spatial curves

\section{Introduction}

A function $F: {\mathbb R}^n\rightarrow {\mathbb R}$ defined on the $n$-dimensional linear space ${\mathbb R}^n$ is called a \emph{convex distance function}, or a \emph{gauge}, if it satisfies the following conditions:

(i) $F(x)\geq 0$ for $x \in {\mathbb R}^n$, and $F(x) = 0 \Leftrightarrow x = 0$,

(ii) $F(\lambda x) = \lambda F(x)$ for $x \in {\mathbb R}^n, \lambda > 0$,

(iii) $F(x+y) \leq F(x)+F(y)$ for $x, y \in {\mathbb R}^n$. \\
Having this notion, the space $({\mathbb R}^n, F)$ is called an $n$-dimensional \emph{generalized Minkowski space} or \emph{gauge space} (cf. \cite{JMR}). It turns out that gauge spaces are characterized by the geometric property to have an arbitrary convex body $B$ as unit ball, where the origin $0$ is from the interior of $B$. (If 
$0$ is even the center of symmetry of $B$, then we have the subcase of normed or 
\emph{Minkowski spaces}, see \cite{T} and \cite{MSW} for basic properties of these spaces.)

Continuing the investigations from \cite{BMSa}, which refer to curves in generalized Minkowski planes in the spirit of the normed subcase developed in \cite{BMSh}, we will study now spatial curves in $3$-dimensional gauge spaces. Our paper is organized as follows. In Section 2 we discuss Frenet-type equations and invariants for curves in $3$-dimensional generalized Minkowski spaces. In Section 3 we characterize generalized cylindrical helices and rectifying curves in such spaces, and in Section 4 we investigate how to compute the invariants and give examples in a Randers space. In Section 5 we discuss how invariants of spatial curves are changing when the unit sphere changes its position with respect to $0$ via parallel translation. It should be noticed that the paper \cite{G} contains related (but mostly different) results; the differences and relations between the results in \cite{G} and our results are always described in detail where they occur.

\section{Frenet-type equations and invariants}

Let $({\mathbb R^3}, F)$ be a $3$-dimensional generalized Minkowski space whose \emph{unit ball} $B$ and whose \emph{unit sphere} $S$ are defined by
\[B = \{x \in {\mathbb R}^3; F(x) \leq 1\}, \ \ \ \ S = \{x \in {\mathbb R}^3; F(x) = 1\}.\]
Here $B$ is a compact convex set having the origin $0$ as interior point and $S$ as its boundary. In the following we assume that $S$ is \emph{smooth} and \emph{strictly convex}, i.e., that through any $x \in S$ a unique supporting plane of $B$ passes, and that $S$ contains no line segment.

Let $v$ be a non-zero vector in ${\mathbb R^3}$, and let $H$ be an oriented plane in ${\mathbb R^3}$ with oriented basis $\{X, Y\}$. Taking the orientation into consideration, we say that $v$ is \emph{Birkhoff orthogonal} to $H$ (denoted by $v \dashv_{B} H$) if the tangent plane of $S$ at $v/F(v)$ is $H$ and $\{v, X, Y\}$ is positively oriented (cf. \cite{BMT} for the subcase of normed spaces).

Let $\gamma(s)$ be an \emph{oriented smooth curve} in $({\mathbb R^3}, F)$ with \emph{arc-length parameter} $s$ such that $F(\gamma'(s)) = 1$. Set $e_1 = \gamma'(s) \in S$, and assume that $e_1' \neq 0$. Let $H$ be the oriented osculating plane spanned by $e_1$ and $e_1'$ such that $\{e_1, e_1'\}$ is an oriented basis. Since $S$ is strictly convex, there exists a unique vector $v = v(s) \in S$ such that $v \dashv_{B} H$.

We can write
\[e_1'' = p_1 e_1+p_2 e_1'+p_3 v\]
for some functions $p_i(s)$. For a positive function $k(s)$,
\[\left( \frac{e_1'}{k} \right)' = \frac{p_1}{k} e_1+\frac{p_2 k-k'}{k^2} e_1'+\frac{p_3}{k} v.\]
If we choose $k(s)$ like
\[k(s) = c_1 \exp\left( \int_{s_{0}}^{s} p_2(\sigma) d\sigma \right)\ \]
where $c_1$ is a positive constant, then
\[\left( \frac{e_1'}{k} \right)' = \frac{p_1}{k} e_1+\frac{p_3}{k} v.\]
Set $e_2 = e_1'/k$. Then
\[e_1' = k e_2, \ \ \ \ e_2' = \frac{p_1}{k} e_1+\frac{p_3}{k} v.\]

By construction, we have $v' \in H$. So $v'$ is spanned by $e_1$ and $e_1'$. We can write
\[v' = q_1 e_1+q_2 e_1' \]
for some functions $q_j(s)$. For a function $f(s)$,
\[(v+f e_1)' = (f'+q_1)e_1+k(f+q_2)e_2.\]
If we choose $f(s)$ by
\[f(s) = -\int_{s_0}^{s} q_1(\sigma)d\sigma+c_2, \]
where $c_2$ is a constant, then
\[(v+f e_1)' = k(f+q_2)e_2.\]
Set $e_3 = v+f e_1$. Then
\[e_3' = k(f+q_2)e_2.\]

Now we have
\[e_1' = k e_2, \ \ \ \ e_2' = \frac{p_1-f p_3}{k} e_1+\frac{p_3}{k} e_3, \ \ \ \ e_3' = k(f+q_2)e_2.\]
Set
\[k^{\ast} = \frac{f p_3-p_1}{k}, \ \ \ \ w = \frac{p_3}{k}, \ \ \ \ w^{\ast} = -k(f+q_2).\]
Then
\[e_1' = k e_2, \ \ \ \ e_2' = -k^{\ast} e_1+w e_3, \ \ \ \ e_3' = -w^{\ast} e_2.\]

\begin{thm}
With the above derivation we have the Frenet-type equations
\[e_1' = k e_2, \ \ \ \ e_2' = -k^{\ast} e_1+w e_3, \ \ \ \ e_3' = -w^{\ast} e_2 \]
for curves in $3$-dimensional generalized Minkowski spaces. There is a $2$-dimensional freedom of $c_1 > 0$ and $c_2$, and the coefficients $k$, $k^{\ast}$, $w$ and $w^{\ast}$ are not invariant.
\end{thm}

{\bf Remark.} For the Frenet-type equations of curves in $({\mathbb R}^3, F)$ studied in \cite{G}, the choice of the frame is restricted in such a way that the total frame spans the unit volume, and so there is a $1$-dimensional freedom.\\

Corresponding to the $2$-dimensional freedom of $c_1 > 0$ and $c_2$, we can write another choice of $k$ and $f$ by
\[\bar{k} = a k, \ \ \ \ \bar{f} = f+b,\]
where $a > 0$ and $b$ are constant. Correspondingly, we obtain
\[\bar{e}_2 = \frac{1}{a} e_2, \ \ \ \ \bar{e}_3 = e_3+b e_1.\]
Computing $\bar{e}_2'$ and $\bar{e}_3'$, we can find that
\[\bar{k}^{\ast} = \frac{1}{a}(k^{\ast}+bw), \ \ \ \ \bar{w} = \frac{1}{a} w, \ \ \ \ \bar{w}^{\ast} = a(w^{\ast}-bk).\]
Thus we get the following theorem (cf. \cite{G}).

\begin{thm}
For a curve $\gamma(s)$ in  $(R^3, F)$ given as above, the quantities 
\[kw, \ \ \ \ \frac{k'}{k}, \ \ \ \ kk^{\ast}+ww^{\ast}, \ \ \ \ \left(\frac{w^{\ast}}{k}\right)' \]
are invariants.
\end{thm}

{\bf Remark.} When $w \neq 0$, then $w'/w$ and $(k^{\ast}/w)'$ are invariants. Since
\[(\log |kw|)' = \frac{k'}{k}+\frac{w'}{w}, \ \ \ \ \left(\frac{kk^{\ast}+ww^{\ast}}{kw} \right)' = \left(\frac{w^{\ast}}{k} \right)'+\left(\frac{k^{\ast}}{w} \right)',\]
$w'/w$ and $(k^{\ast}/w)'$ are given by combinations of the four invariants in Theorem 2.2.\\

{\bf Remark.} From the construction, we can see that
\[kw = p_3, \ \ \ \ \frac{k'}{k} = p_2, \ \ \ \ kk^{\ast}+ww^{\ast} = -p_1-p_3 q_2, \]
\[\left(\frac{w^{\ast}}{k}\right)' = q_1-q_2'.\]

\vspace{3mm}

For convenience, we use the following notation of the invariants:
\[I_1 = kw, \ \ \ \ I_2 = \frac{k'}{k}, \ \ \ \ I_3 = kk^{\ast}+ww^{\ast}, \ \ \ \ I_4 = \left(\frac{w^{\ast}}{k}\right)'. \]

\section{Cylindrical helices and rectifying curves}

It is well known that, by several reasons, cylindrical helices play an important role in the classical differential geometry of spatial curves.: In Euclidean $3$-space, a curve $\gamma$ is called a \emph{cylindrical helix} if there exists a fixed direction such that the angle between all tangential lines of $\gamma$ and this fixed direction is constant, and such a curve is characterized by the condition that the ratio of curvature and torsion is constant. In $({\mathbb R^3}, F)$, the notion of angle is not given; compensating this it makes sense to consider the following framework.

The plane spanned by $\{e_1(s), e_3(s)\}$ is invariant, and we call it the \emph{rectifying plane at} $\gamma(s)$. Using this notion, we get the following theorem, which can be seen as \emph{characterization of cylindrical helices in generalized Minkowski spaces}.

\begin{thm}
There exists a line $\ell$ such that every rectifying plane of a curve $\gamma(s)$ in $({\mathbb R^3}, F)$ contains a line parallel to $\ell$ if and only if $(w^{\ast}/k)' = 0$.
\end{thm}

Proof. For our proof we fix a moving frame.

(i) Suppose that there exists a line $\ell$ such that every rectifying plane contains a line parallel to $\ell$. We fix a non-zero vector $a_0$ which is parallel to $\ell$. Then there exist functions $\alpha(s)$ and $\beta(s)$ such that
\[\alpha e_1+\beta e_3 = a_0.\]
Differentiating, we get
\[\alpha' e_1+(\alpha k-\beta w^{\ast})e_2+\beta' e_3 = 0.\]
So $\alpha$ and $\beta$ are constant and $\alpha k = \beta w^{\ast}$. Since $k > 0$, we have that if $\beta = 0$, then $\alpha = 0$, which is a contradiction. So $\beta \neq 0$ and
\[\frac{w^{\ast}}{k} = \frac{\alpha}{\beta}, \ \ \ \ \left( \frac{w^{\ast}}{k} \right)' =0.\]

(ii) Suppose that $(w^{\ast}/k)' = 0$. Then $w^{\ast}/k$ is constant, which we set $c$. Since
\[(ce_1+e_3)' = (ck-w^{\ast})e_2 = 0 \]
and $\{e_1, e_3\}$ are linearly independent, we have $ce_1+e_3 = a_0$ for some fixed non-zero vector $a_0$. Thus the converse is verified, too. \ \ \ \ $\Box$ \\

A curve $\gamma(s)$ in $({\mathbb R^3}, F)$ is called a \emph{rectifying curve} if there exists a fixed point $p_0$ such that every rectifying plane contains $p_0$ (cf. \cite{C} for the Euclidean case). We extend the characterization of rectifying curves given in \cite{C} to (norms and) gauges as follows.

\begin{thm}
A curve $\gamma(s)$ in $({\mathbb R^3}, F)$ is a rectifying curve if and only if $(w^{\ast}/k)'$ is a non-zero constant.
\end{thm}

Proof. (i) Suppose that $\gamma(s)$ is a rectifying curve. Then there exists a fixed point $p_0$ such that
\[\gamma+\alpha e_1+\beta e_3 = p_0 \]
for some functions $\alpha(s)$ and $\beta(s)$. Differentiating, we obtain
\[(1+\alpha')e_1+(\alpha k-\beta w^{\ast})e_2+\beta' e_3 = 0.\]
So $\beta$ is constant, $\alpha = c-s$ for a constant $c$, and $(c-s)k = \beta w^{\ast}$. Noting that $\beta \neq 0$, we have
\[\frac{w^{\ast}}{k} = \frac{c-s}{\beta}, \ \ \ \ \left( \frac{w^{\ast}}{k} \right)' = -\frac{1}{\beta}.\]

(ii) Suppose that $(w^{\ast}/k)'$ is a non-zero constant, which we set $-1/\beta$. Then
\[\frac{w^{\ast}}{k} = \frac{c-s}{\beta} \]
for a constant $c$. Since
\[(\gamma+(c-s)e_1+\beta e_3)' = \{ (c-s)k-\beta w^{\ast} \}e_2 = 0,\]
we have $\gamma+(c-s)e_1+\beta e_3 = p_0$ for some fixed point $p_0$. Thus, also the converse is proved. \ \ \ \ $\Box$

\section{Examples in a Randers space}

Let $\gamma(t) = (\gamma_1(t), \gamma_2(t), \gamma_3(t))$ be an oriented smooth curve in $({\mathbb R}^3, F)$ with arbitrary parameter $t$, and assume that $\{\gamma'(t), \gamma''(t)\}$ are linearly independent. The oriented osculating plane $H$ satisfies $H = \mbox{span}\{\gamma'(t), \gamma''(t)\}$ with oriented basis $\{\gamma'(t), \gamma''(t) \}$. We denote by $\times$ the standard cross product in ${\mathbb R}^3$. The vector $v \in S$ with $v \dashv_{B} H$ is characterized by the condition
\[\mbox{grad}(F)|_v = \lambda \gamma' \times \gamma'', \ \ \ \ \lambda > 0.\]
But it is not easy to treat this in general.

So we shall consider the case of Randers spaces (cf. \cite{BCS} and \cite{CS}). Let $F$ be a \emph{Randers norm} given by
\[F(x_1, x_2, x_3) = \sqrt{x_1^2+x_2^2+x_3^2}+bx_3, \ \ \ \ 0 < b < 1.\]
Then the unit sphere $S$ is a rotational ellipsoid given by
\[G(x_1, x_2, x_3):= (1-b^2)(x_1^2+x_2^2)+\{(1-b^2)x_3+b\}^2 = 1.\]
Using $G$ instead of $F$, we see that the vector $v \in S$ with $v \dashv_{B} H$ is characterized by the condition
\[\mbox{grad}(G)|_v = \lambda \gamma' \times \gamma'', \ \ \ \ \lambda > 0.\]
Let
\[\gamma' \times \gamma'' =: (y_1, y_2, y_3).\]
Then the vector $v \in S$ with $v \dashv_{B} H$ is given by
\[v = \left( \begin{array} {c} 
\frac{y_1}{\sqrt{(1-b^2)(y_1^2+y_2^2)+y_3^2}} \\
\frac{y_2}{\sqrt{(1-b^2)(y_1^2+y_2^2)+y_3^2}}\\
\frac{1}{1-b^2}\left(\frac{y_3}{\sqrt{(1-b^2)(y_1^2+y_2^2)+y_3^2}}-b \right)
\end{array} \right).\]
Let $s$ be the arc-length parameter of $\gamma$. Then
\[e_1 = \frac{d\gamma}{ds} = \frac{\gamma'}{F(\gamma')},\]
and we can write
\[\frac{d^2 e_1}{ds^2} = p_1 e_1+p_2 \frac{de_1}{ds}+p_3 v,\]
\[\frac{dv}{ds} = q_1 e_1+q_2 \frac{de_1}{ds}. \]
To compute the invariants, it is necessary to calculate the coefficients.\\

{\bf Example 1.} In the above Randers space $({\mathbb R}^3, F)$, let
\[\gamma(t) = \frac{1}{\sqrt{2}+b}(\cos t, \sin t, t).\]
Then
\[\gamma' = \frac{1}{\sqrt{2}+b}(-\sin t, \cos t, 1).\]
Since $F(\gamma') = 1$, we find that $t$ is the arc-length parameter and $e_1 = \gamma'$. We compute
\[e_1' = \gamma'' = \frac{1}{\sqrt{2}+b}(-\cos t, -\sin t, 0),\]
\[e_1'' = \frac{1}{\sqrt{2}+b}(\sin t, -\cos t, 0), \]
\[\gamma' \times \gamma'' = \frac{1}{(\sqrt{2}+b)^2} (\sin t, -\cos t, 1),\]
and
\[v = \left( \frac{\sin t}{\sqrt{2-b^2}}, \ -\frac{\cos t}{\sqrt{2-b^2}}, \ \frac{1}{1-b^2}\left(\frac{1}{\sqrt{2-b^2}}-b\right) \right).\]
Hence we get
\[v' = -\frac{\sqrt{2}+b}{\sqrt{2-b^2}} e_1',\]
and so
\[q_1 = 0, \ \ \ \ q_2 =  -\frac{\sqrt{2}+b}{\sqrt{2-b^2}}.\]
Thus we have $I_4 = 0$. In fact, every rectifying plane contains a line which is parallel to a fixed vector $(0,0,1)$.

We can also compute
\[e_1'' = -\frac{\sqrt{2-b^2}-b}{2\sqrt{2-b^2}} e_1+\frac{\sqrt{2-b^2}+b}{2(\sqrt{2}+b)} v, \]
and so
\[p_1 = -\frac{\sqrt{2-b^2}-b}{2\sqrt{2-b^2}}, \ \ \ \ p_2 = 0,\ \ \ \ p_3 = \frac{\sqrt{2-b^2}+b}{2(\sqrt{2}+b)}.\]
Thus we have
\[I_1 = \frac{\sqrt{2-b^2}+b}{2(\sqrt{2}+b)}, \ \ \ \ I_2 = 0, \ \ \ \ I_3 = 1.\]
The four invariants are all constant in this case.\\

{\bf Example 2.} Let $a(t) = (a_1(t), a_2(t), a_3(t))$ be an oriented smooth curve in the above Randers space $({\mathbb R}^3, F)$ with arc-length parameter $t$, and assume that $a''(t) \neq 0$. For a positive function $f(t)$, let
\[\gamma(t) = \int f(t)a'(t)dt = \left( \int f(t) a_1'(t)dt , \int f(t) a_2'(t)dt , \int f(t) a_3'(t)dt \right).\]
Let $s$ be the arc-length parameter of $\gamma$. Then 
\[\gamma' = fa', \ \ \ \ F(\gamma') = f = \frac{ds}{dt}, \ \ \ \ e_1 = \frac{d\gamma}{ds} = a',\]
\[\frac{de_1}{ds} = \frac{1}{f}a'', \ \ \ \ \frac{d^2 e_1}{ds^2} = -\frac{f'}{f^3}a''+\frac{1}{f^2}a^{(3)},\]
\[\gamma'' = f'a'+fa'', \ \ \ \ \gamma' \times \gamma'' = f^2 a'\times a''.\]
The oriented osculating plane of $\gamma$ is the same as that of $a$. So the vector $v$ for $\gamma$ is the same as that for $a$, which is independent of $f$. We can write
\[a^{(3)} = P_1 a'+P_2 a''+P_3 v, \ \ \ \ v' = Q_1 a'+Q_2 a'', \]
where $P_i(s)$ and $Q_j(s)$ are functions which are independent of $f$. Let $I^a_i$ $(1 \leq i \leq 4)$ denote the invariants of $a$. Then
\[I^a_1 = P_3, \ \ \ \ I^a_2 = P_2, \ \ \ \ I^a_3 = -P_1-P_3 Q_2, \ \ \ \ I^a_4 = Q_1-Q_2'.\]
We have 
\[\frac{d^2 e_1}{ds^2} = \frac{P_1}{f^2} e_1+\left( \frac{P_2}{f}-\frac{f'}{f^2} \right) \frac{de_1}{ds}+\frac{P_3}{f^2}v,\]
and
\[\frac{dv}{ds} = \frac{Q_1}{f} e_1+Q_2 \frac{de_1}{ds}.\]
Thus, for $\gamma$ we get
\[p_1 = \frac{P_1}{f^2}, \ \ \ \ p_2 = \frac{P_2}{f}-\frac{f'}{f^2}, \ \ \ \ p_3 = \frac{P_3}{f^2}, \ \ \ \ q_1 = \frac{Q_1}{f}, \ \ \ \ q_2 = Q_2,\]
and the invariants are as follows:
\[I_1 = \frac{P_3}{f^2} = \frac{I^a_1}{f^2}, \ \ \ \ I_2 = \frac{P_2}{f}-\frac{f'}{f^2} = \frac{I^a_2}{f}-\frac{f'}{f^2},\]
\[I_3 = -\frac{P_1+P_3 Q_2}{f^2} = \frac{I^a_3}{f^2},\ \ \ \ I_4 = q_1-\frac{dq_2}{ds} = \frac{Q_1-Q_2'}{f} = \frac{I^a_4}{f}. \]
From this we can see the following

\begin{prop}
In Example 2, the following statements hold:

(i) If $\{a', a'', a^{(3)}\}$ are linearly independent ($P_3 \neq 0$), then there exists a function $f$ such that $I_1$ is a non-zero constant.

(ii) There exists a function $f$ such that $I_2$ is constant.

(iii) If $I^a_3 \neq 0$, then there exists a function $f$ such that $I_3$ is a non-zero constant.

(iv) If $I^a_4 \neq 0$, then there exists a function $f$ such that $I_4$ is a non-zero constant.
\end{prop}

\section{Gauge change by parallel translations}

As in Section 2, let $\gamma(s)$ be an oriented smooth curve in $({\mathbb R^3}, F)$ with arc-length parameter $s$ such that $F(\gamma'(s)) = 1$. Set $e_1 = \gamma'(s) \in S$, and assume that $e_1' \neq 0$. Let $H$ be the oriented osculating plane, and there exists a unique vector $v = v(s) \in S$ such that $v \dashv_{B} H$.

Let $\bar{F}$ be another gauge whose unit sphere $\bar{S}$ satisfies $\bar{S} = S+a_0$ for a fixed non-zero vector $a_0$. Let $\bar{s}$ be the arc-length parameter of $\gamma$ with respect to $\bar{F}$. Then
\[\bar{e}_1 = \frac{d\gamma}{d\bar{s}} = \frac{d\gamma}{ds} \ \frac{ds}{d\bar{s}} = f e_1,\]
where we set $f = ds/d{\bar{s}} > 0$. The oriented osculating plane $H$ is not changed, and the vector $\bar{v} \in \bar{S}$ which is Birkhoff orthogonal to $H$ with respect to $\bar{F}$ satisfies $\bar{v} = v+a_0$. 

Then
\[\frac{d\bar{e}_1}{d\bar{s}} = \frac{d\bar{e}_1}{ds} \ \frac{ds}{d\bar{s}} = f \left(f' e_1+f e_1' \right)\]
and
\[\frac{d\bar{v}}{d\bar{s}} = \frac{dv}{d\bar{s}} = \frac{dv}{ds} \ \frac{ds}{d\bar{s} }= f (q_1 e_1+q_2 e_1') = \left(q_1-\frac{q_2 f'}{f} \right)\bar{e}_1+\frac{q_2}{f} \ \frac{d\bar{e}_1}{d\bar{s}}.\]
So
\[\bar{q}_1 = q_1-\frac{q_2 f'}{f}, \ \ \ \ \bar{q}_2 = \frac{q_2}{f},\]
and
\[\bar{I}_4 = \bar{q}_1-\frac{d\bar{q}_2}{d\bar{s}} = q_1-\frac{q_2 f'}{f}-f\left(\frac{q_2}{f}\right)' = q_1-q_2' = I_4.\]
Thus we get the following

\begin{thm}
The quantity $I_4$ is not changed under parallel translations of unit spheres.
\end{thm}

For the other quantities $I_1, I_2$ and $I_3$ we compute
\[\frac{d^2 \bar{e}_1}{d\bar{s}^2} = (ff''-2f'^2+p_1 f^2-p_2 ff')\bar{e}_1+(3f'+p_2 f)\frac{d\bar{e}_1}{d\bar{s}}+p_3 f^3(\bar{v}-a_0).\]
If we write
\[a_0 = a_{01}\bar{e}_1+a_{02}\frac{d\bar{e}_1}{d\bar{s}}+a_{03}\bar{v},\]
then
\[\bar{p}_1 = ff''-2f'^2+p_1 f^2-p_2 ff'-p_3 f^3 a_{01},\]
\[\bar{p}_2 = 3f'+p_2 f-p_3 f^3 a_{02}, \ \ \ \ \bar{p}_3 = p_3 f^3(1-a_{03}).\]
Using these translation formulas, we will discuss examples.

Let $F$ be a Randers norm given by
\[F(x_1, x_2, x_3) = \sqrt{x_1^2+x_2^2+x_3^2}+bx_3, \ \ \ \ 0 < b < 1.\]
Then the unit sphere $S$ is
\[(1-b^2)(x_1^2+x_2^2)+\{(1-b^2)x_3+b\}^2 = 1.\]
If
\[a_0 = \left(0, 0, \frac{b}{1-b^2}\right),\]
then $\bar{S} = S+a_0$ is given by
\[(1-b^2)(x_1^2+x_2^2)+(1-b^2)^2 x_3^2 = 1,\]
and, correspondingly,
\[\bar{F}(x_1, x_2, x_3) = \sqrt{(1-b^2)(x_1^2+x_2^2)+(1-b^2)^2 x_3^2}.\]
In the following examples, we use these gauges $F$ and $\bar{F}$.\\

{\bf Example 3.} As in Example 1 in Section 4, let
\[\gamma(t) = \frac{1}{\sqrt{2}+b}(\cos t, \sin t, t).\]
Then
\[\gamma' = \frac{1}{\sqrt{2}+b}(-\sin t, \cos t, 1)\]
and $F(\gamma') = 1$, such that $t$ is the arc-length parameter with respect to $F$. As shown in Example 1, we have
\[v = \left( \frac{\sin t}{\sqrt{2-b^2}}, \ -\frac{\cos t}{\sqrt{2-b^2}}, \ \frac{1}{1-b^2}\left(\frac{1}{\sqrt{2-b^2}}-b\right) \right), \]
\[q_1 = 0, \ \ \ \ q_2 = -\frac{\sqrt{2}+b}{\sqrt{2-b^2}},\]
\[p_1 = -\frac{\sqrt{2-b^2}-b}{2\sqrt{2-b^2}}, \ \ \ \ p_2 = 0,\ \ \ \ p_3 = \frac{\sqrt{2-b^2}+b}{2(\sqrt{2}+b)},\]
and
\[I_1 = \frac{\sqrt{2-b^2}+b}{2(\sqrt{2}+b)}, \ \ \ \ I_2 = 0, \ \ \ \ I_3 = 1.\]

With respect to $\bar{F}$, we have
\[\bar{F}(\gamma') = \frac{\sqrt{(1-b^2)(2-b^2)}}{\sqrt{2}+b} = \frac{d\bar{s}}{dt},\]
where $\bar{s}$ is the arc-length parameter for $\bar{F}$. Then
\[\frac{dt}{d\bar{s}} = \frac{\sqrt{2}+b}{\sqrt{(1-b^2)(2-b^2)}} =: f,\]
\[\bar{e}_1 = \frac{1}{\sqrt{(1-b^2)(2-b^2)}} (-\sin t, \cos t, 1),\]
and
\[\bar{v} = v+a_0 = \frac{1}{\sqrt{2-b^2}}\left(\sin t, -\cos t, \frac{1}{1-b^2} \right).\]
Since
\[a_0 = b\sqrt{\frac{1-b^2}{2-b^2}}\bar{e}_1+\frac{b}{\sqrt{2-b^2}}\bar{v},\]
we have
\[a_{01} = b\sqrt{\frac{1-b^2}{2-b^2}}, \ \ \ \ a_{02} = 0, \ \ \ \ a_{03} = \frac{b}{\sqrt{2-b^2}}.\]
Hence we obtain
\[\bar{q}_1 = 0, \ \ \ \ \bar{q}_2 = -\sqrt{1-b^2},\]
\[\bar{p}_1 = -\frac{(\sqrt{2}+b)^2}{(1-b^2)(2-b^2)^2},\ \ \ \ \bar{p}_2 = 0, \ \ \ \ \bar{p}_3 = \frac{(\sqrt{2}+b)^2}{\sqrt{1-b^2}\ (2-b^2)^2},\]
and
\[\bar{I_1} = \frac{(\sqrt{2}+b)^2}{\sqrt{1-b^2}\ (2-b^2)^2}, \ \ \ \ \bar{I}_2 = 0, \ \ \ \ \bar{I}_3 = \frac{(\sqrt{2}+b)^2}{(1-b^2)(2-b^2)}. \]

So the quantities $I_1$ and $I_3$ change under parallel translations of unit spheres. Also by a curve which is $C^3$-close to $\gamma$, we can see that $I_1$ and $I_3$ change under parallel translations of unit spheres.\\

{\bf Example 4.} Let
\[\gamma(t) = \frac{1}{1-b^2}(0, \sqrt{1-b^2}\cos t, \sin t).\]
Then
\[\gamma' =  \frac{1}{1-b^2}(0, -\sqrt{1-b^2}\sin t, \cos t)\]
and $\bar{F}(\gamma') = 1$. So $t$ is the arc-length parameter $\bar{s}$ with respect to $\bar{F}$, and $\bar{e}_1 = \gamma'$. Since $\bar{e}_1'' = -\bar{e}_1$, we have $\bar{p}_2 = 0 = \bar{I}_2$.

With respect to $F$, we have
\[F(\gamma') = \frac{1}{1-b^2}(\sqrt{1-b^2 \sin^2{t}}+b\cos t) = \frac{ds}{dt},\]
where $s$ is the arc-length parameter of $\gamma$ for $F$. Set
\[g := \frac{dt}{ds} = \frac{1-b^2}{\sqrt{1-b^2 \sin^2{t}}+b\cos t} = \sqrt{1-b^2 \sin^2{t}}-b\cos t.\]
Then
\[e_1 = g\gamma', \ \ \ \ \frac{de_1}{ds} = gg'\gamma'+g^2\gamma'',\]
\[\frac{d^2 e_1}{ds^2} = g(g'^2+gg''-g^2)\gamma'+3g^2 g'\gamma''\]
\[= (gg''-g^2-2g'^2)e_1+3g'\frac{de_1}{ds},\]
and we can see that $p_2 = 3g' = I_2$.

So the quantity $I_2$ is changed under parallel translations of unit spheres. Also by a curve which is $C^3$-close to $\gamma$, we can see that $I_2$ is changed under parallel translations of unit spheres.\\

Theorem 5.1 as well as Examples 3 and 4 imply that the quantity $I_4$ is an invariant of higher degree than $I_1, I_2$ and $I_3$.

\end{document}